\documentclass[12pt, reqno]{amsart}
\usepackage{amsmath}
\usepackage{amsfonts}
\usepackage{amssymb}
\usepackage[english]{babel}
\theoremstyle{plain}
\newtheorem {theorem}{Theorem}[section]

\newtheorem {proposition} [theorem]{Proposition}
\theoremstyle{definition}

\newtheorem{remark}[theorem]{Remark}
\newtheorem{example}[theorem]{Example}
\theoremstyle{remark}

\renewcommand{\rho}{\varrho}
\renewcommand{\epsilon}{\varepsilon}
\renewcommand{\phi}{\varphi}
\renewcommand{\theta}{\vartheta}

\newcommand{\la}{\lambda}
\newcommand{\ga}{\gamma}

\renewcommand{\d}{\delta}

\renewcommand{\a}{\alpha}
\renewcommand{\b}{\beta}

\newcommand{\g}{\gamma}

\newcommand{\p}{\partial}

\renewcommand{\gg}{{\bar\g}}

\newcommand{\barga}{{\bar\ga}}

\newcommand{\bb}{{\bar\b}}
\renewcommand{\aa}{{\bar\a}}

\newcommand{\hh}{{\bar{h}}}
\newcommand{\kk}{{\bar{k}}}

 \newcommand{\mm}{{\bar\mu}}
\newcommand{\barb}{{\bar\b}}

\renewcommand{\ll}{{\bar\la}}

\newcommand{\R}{{\mathbb{R}}}
\newcommand{\C}{{\mathbb{C}}}

\newcommand{\G}{{\Gamma}}

\newcommand{\NN}{{\bar{N}}}

\renewcommand{\H}{{\mathcal{H}}}
\newcommand{\HH}{{{\bar{\mathcal{H}}}}}
\newcommand{\VV}{{\bar{V}}}
\newcommand{\UU}{{\bar{U}}}
\newcommand{\WW}{{\bar{W}}}
\newcommand{\ZZ}{{\bar{Z}}}

\newcommand{\D}{{\nabla}}

\newcommand{\barint}
{\rule[.036in]{.12in}{.009in}\kern-.16in \displaystyle\int}

\title[Levi umbilical surfaces in complex space  ]{ Levi umbilical surfaces
in complex space}

\author{Roberto Monti}
\address{Roberto Monti: Dipartimento di Matematica Pura ed Applicata \\
Universit\`a di Padova \\ Via Belzoni, 7 \\35100 Padova, Italy}
 \email{%rmonti@science.unitn.it,
monti@math.unipd.it}

\author{Daniele Morbidelli}
\address{Daniele Morbidelli: Dipartimento di Matematica\\
 Universit\`a di Bologna\\
 Piazza di Porta San Donato, 5\\
 40127  Bologna, Italy}
\email{morbidel@dm.unibo.it}

\date{1st November 2005}

\begin{document}

\begin{abstract}
We define a complex connection on a real hypersurface of
$\C^{n+1}$ which is naturally inherited from the ambient space.
Using a system of Codazzi-type equations, we classify connected
real hypersurfaces in $\C^{n+1}$, $n\ge 2$, which are Levi
umbilical and have non zero constant Levi curvature. It turns out
that such surfaces are contained either in a sphere or in the
boundary of a complex tube domain with spherical section.
\end{abstract}

\maketitle

\section{Introduction}

Let $M$ be a $(2n+1)$-dimensional real surface embedded in
$\C^{n+1}$, denote by $h$ the $\C$-linear extension of the second
fundamental form of $M$ and by  $g$ be the restriction to the
complexified tangent bundle $\C T M$  of the standard hermitian
product of $\C^{n+1}$. The surface $M$ is Levi umbilical  if $h(Z,
\bar W)= H g(Z, \bar W)$ for some scalar function $H$ (the Levi
curvature) and for all holomorphic  tangent vector fields $Z$ and
$W$. Levi umbilicality is weaker than Euclidean umbilicality
because it contains no information on terms of the form $h(Z,W)$
with holomorphic $Z$ and $W$. In particular, it is easy to
construct Levi umbilical surfaces which are neither spheres nor
hyperplanes. Indeed, any surface which is the zero set $F=0$ of a
smooth defining function $F(z,\bar z)=|z|^2+\Phi(z,\bar z)$, where
$\Phi$ is any polyharmonic function in $\C^{n+1}$, is Levi
umbilical (see Example \ref{es2}).

In view of these examples, a natural question is whether there is
any version of the classical  Darboux theorem for usual umbilical
surfaces. In this paper  we classify Levi umbilical surfaces with
\emph{constant} non zero Levi curvature. An example of such
surfaces are, of course, the spheres $\{z\in\C^{n+1}:|z|=r\}$,
$r>0$. A less trivial example is the boundary of spherical tubes,
i.e.~surfaces of the form (see Example \ref{esso})
\begin{equation}
 \label{cyl}
      \Big\{ z\in\C^{n+1}: \sum_{h=1}^{n+1} (z_h +\bar z_h)^2=
      r^2\Big\},\quad r>0.
\end{equation}

\noindent Our main result states that there are no other examples.
More precisely, we prove that any $(2n+1)$-dimensional oriented
connected surface embedded in $\C^{n+1}$, $n\geq 2$, which is Levi
umbilical and has non zero constant Levi curvature  is necessarily
contained either in a sphere or, up to complex isometries of
$\C^{n+1}$, in a spherical cylinder of the form \eqref{cyl}. This
is proved in Theorem \ref{teor}. It is interesting to observe the
appearance of tube domains, which are relevant objects in several
complex variables, see \cite{Kr}.

This classification follows from the analysis of a system of
Codazzi  equations for $h$, where covariant derivatives are
computed with respect to a suitable complex connection $\D$ on
$M$. Though very natural, this connection and the corresponding
Codazzi   equations do not seem to be studied in the literature.
The main features of $\D$ are:
\begin{itemize}
\item
[(a)] both  the holomorphic and the antiholomorphic bundles are
parallel;

\item
[(b)] the restriction  $g$ to $\C TM$ of the hermitian product in
$\C^{n+1}$ satisfies $\D g=0$.
\end{itemize}

\noindent Briefly, the connection is constructed in the following
way. Let $\nu$ be a real unit normal to $M$ and consider $N =
2^{-1/2}\left( \nu -i T\right)$, the holomorphic unit normal to
$M$. Here, $T= J(\nu)$ where $J$ is the standard complex structure
of $\C^{n+1}$. Then, given a holomorphic tangent vector field $Z$
and a tangent vector $U$, we define
\[
       \D_U Z = D_U Z - g(  D_U Z, \NN)  N,
\]
where $D$ is the standard connection in $\C^{n+1}$. Then, this
definition, along with   $\D T =0$, is extended to the whole
tangent bundle, giving rise to a connection satisfying (a) and (b)
(see Section \ref{cnx}).

Properties (a) and (b) are similar to the ones of the
Tanaka-Webster connection on strictly pseudoconvex Cauchy-Riemann
manifolds (see \cite{T} and \cite{W}). Whereas for this connection
 the Levi form $-i d\theta$ associated with a
contact form $\theta$ plays the role of the metric and is required
to be parallel, in our case the metric inherited from $\C^{n+1}$
is required to be parallel. See also the discussion in
 Remark \ref{webb}. This produces a connection
which seems to be  more suitable for our purposes.  A different
connection is introduced by Klingenberg in \cite{K}. It
 arises as orthogonal projection of the standard
connection in the space and, in general, does not satisfy property
(a).

A typical example of Codazzi equation for $h$, written in
components with respect to a holomorphic frame $Z_1,...,Z_n$, is
(see Remark \ref{remi})
\begin{equation}
 \label{cech}
   \nabla_\a h_{\b\barga} -  \nabla_\b h_{\a\barga}   =
   i   h _{\b\barga}h_{\a0}-  i h _{\a\barga} h_{\b0},
\end{equation}
where $h_{\a\bb}=h(Z_\a, \bar Z_\b)$ for $\a,\b=1,...,n$ and index
$0$ refers to $T$. In Theorem \ref{coda}, we compute the system of
equations needed in the classification theorem. In these
equations, as in \eqref{cech}, there is a non vanishing right hand
side, reflecting both the non vanishing of $\mathrm{Tor}_\D$ and
the non vanishing of $g(D_Z N,\NN)$.

Concerning the restriction $n\geq 2$ in the classification
theorem, note that for $n=1$ the umbilicality property is
satisfied by any hypersurface of $\C^2$.  Moreover,  by the
existence and regularity results proved by Slodkowski and
Tomassini \cite{ST} and Citti, Lanconelli and Montanari \cite{CLM}
for the Levi equation, there are smooth graphs in $\C^2$ with
prescribed boundary and with constant Levi curvature which do not
belong to the classes described above. Then, the natural question
is whether a \emph{compact} surface in $\C^2$ having constant Levi
curvature is necessarily a sphere. This question has been recently
addressed in \cite{HL} by Hounie and Lanconelli, who give an
affirmative answer in the class of Reinhardt domains.

Another result implied by our Codazzi equations is the
classification of connected pseudoconvex surfaces with non zero
constant Levi curvature and vanishing $h_{\a\b}$ (the symmetric
part  of the second fundamental form). Up to complex isometry,
such surfaces are contained in a sphere or in a spherical cylinder
of the form
\[
       \Big\{ z\in\C^{n+1} : \sum_{i=m}^{n+1} |z_i|^2
       =r^2\Big\},\quad r>0, \quad 1\leq m\leq n.
\]
This is established in Theorem \ref{klino}, which improves Theorem
5.2 in \cite{K}, where the result is proved by a global argument
under  compactness and strict pseudoconvexity assumptions (see
Remark \ref{ria}).

The notion of Levi curvature was introduced by Bedford and Gaveau
in \cite{BG} and it has been recently generalized by Montanari and
Lanconelli in \cite{ML}. There is an increasing interest on
problems concerning this curvature, mainly from the point of view
of partial differential equations. Other significant references
are Citti and Montanari \cite{CM}, Huisken and Klingenberg
\cite{HK} and Montanari and Lascialfari \cite{MLa}. The tools
developed in this work could be useful in the study also of other
problems concerning real hypersurfaces in complex space.

Concerning terminology, we  call ``Levi form" the her\-mitian map
$(Z,W) \mapsto h(Z, \WW)$,  with holomorphic $Z$ and $W$. This is
justified by the fact that $h(Z, \WW)$ coincides with the Levi
form  associated with a natural pseudohermitian structure (see
\cite{W} or \cite{JL} for this notion) inherited by $M$ from the
ambient (see the discussion in Section \ref{suo}).

\medskip \noindent\bf Notation. \rm Greek indices $\a,\b$ etc.~run
from $1$ to $n$, Latin indices $h,k$ run from $1$ to $n+1$. We let
$ \p_h = \frac{\p}{\p z_h}$, $\p_{\bar h} = \frac{\p}{\p\bar z_h}$
and $F_h = \p_h F$. $J$ is the standard complex structure and $D$
is the usual connection in $\C^{n+1}$. The standard hermitian
product $g$ in $\C^{n+1}$ is normalized by $ g(\p_h,\p_\kk )
=g(\p_\hh,\p_k )= \frac 12\d_{hk}$, $g(\p_h,\p_k ) =
g(\p_\hh,\p_\kk )=0$, where $\d_{hk}$ is the Kronecker symbol. The
metric tensors $g_{\a\bb}$ and $g^{\a\bb}$, which are related by
$g^{\a\bb}g_{\ga\bb} =\d_{\a\ga}$,  are used to lower and raise
indices, e.g.~${h_{\a}}^\b = g^{\b\gg} h_{\a\gg}$. If $h$ is
symmetric, we equivalently write $h_{\a}^\b={h_{\a}}^\b$. We adopt
the summation convention. If $E$ is a bundle we denote by
$\Gamma(E)$ the sections of $E$. Finally, $[U,V]$ denotes the Lie
bracket of vector fields and $\mathrm{Tor}_\D (U,V)= \D_U V - \D_V
U -[U,V]$ is the torsion of the connection $\D$.

\medskip

\noindent\bf Acknowledgments. \rm We are indebted to
 Ermanno Lanconelli and
Annamaria Montanari for several fruitful conversations on Levi
curvature.

\section{Levi form and examples}
\label{suo}\setcounter{equation}{0}

Let $M\subset\C^{n+1}$ be a real hypersurface oriented by a real
unit normal $\nu$. We denote by $\H = T^{1,0}M$ (resp.~$\bar \H =
T^{0,1}M$) the holomorphic (resp.~antiholomorphic) tangent bundle
of $M$. We restrict the complex structure $J$ to $\H\oplus\bar\H$
and the metric $g$ to $\C T M $. The vector field  $T=J(\nu)$ is
tangent to $M$. Then, the complexified tangent bundle $\C TM$ can
be decomposed as a direct sum $\H \oplus \bar \H  \oplus \C T$ and
the decomposition is orthogonal with respect to $g$. The
\emph{holomorphic unit normal} to $M$ is the holomorphic vector
field
\begin{equation}
 \label{enne}
              N=\frac{1}{\sqrt 2}\left(\nu-iT\right).
\end{equation}
Up to orientation, $N$ is defined uniquely on $M$ by $|N|=1$ and
$g( N, U)=0$ for all $U\in \H\oplus\bar\H$. Here and in the
following, $|V|^2=g(V,\bar V)$. We have the relations
\begin{equation}
 \label{NT}
  T=\frac{i}{\sqrt 2}(N-\bar N), \quad
  \nu = \frac{1}{\sqrt 2}(N+\bar N).
\end{equation}
There is a unique real $1$-form $\eta$ on $M$ such that
\begin{equation}
\label{pupo}
    \eta(T)=1   \quad \textrm{and} \quad
    \eta(Z)= 0 \quad   \text{for all $Z\in \H\oplus\bar\H.$}
\end{equation}
Precisely, $\eta(Z) = g( Z, T)$ for any $Z\in\C TM$. The
\emph{Levi form} on $M$ associated with $\eta$ is the hermitian
form on $\H$ defined by
\begin{equation}
 \label{Levi}
   L_\eta(Z,\bar W)=\frac {1}{2i}  d\eta (Z, \bar W),
          \quad  Z,W\in \H.
\end{equation}

Denote by $h$ the $\C$-linear extension to $\C TM\times \C TM$ of
the second fundamental form of $M$.  For $Z,W\in \C TM$ let
\begin{equation}
 \label{dh}
      h(Z,W)  = g ( Z, D_{W}\nu) .
\end{equation}
Note that $h(Z,W)=h(W,Z)$ and $\overline{h(Z,W})= h(\bar Z, \bar
W)$.

The Levi form associated with $\eta$ coincides with the hermitian
part of the second fundamental form, i.e.~$L_\eta(Z, \bar W) =  h
(Z, \WW)$ for all $Z,W\in \H$. Indeed, by \eqref{pupo} and
\eqref{NT},
\begin{equation}
 \label{campanella}
\begin{split}
  d\eta (Z,\WW)&= Z \eta(\WW)-\WW\eta(Z) -\eta ([Z,\WW])
                = -\eta([Z,\WW]) \\ &
                =  g ([\WW,Z],T)
                =  \frac{i}{\sqrt 2 } g(D_\WW Z- D_Z\WW,N-\NN).
\end{split}
\end{equation}
Since $g(D_\WW Z,N)=g(D_Z\WW,\NN)=0$, we find
\begin{equation}
 \label{campanello}
   d\eta (Z,\WW) =- \frac{i}{\sqrt 2 }
                     \big( g(D_\WW Z,N+\NN)
                        +g(D_Z\WW,N+\NN) \big)
          =- 2i g(D_\WW Z,\nu).
\end{equation}
The claim follows.

The \emph{Levi curvature} $H$ of $M$ is the trace of the Levi
form. The surface $M$ is \emph{Levi umbilical} if $h (Z,\WW)= H
g(Z,\WW)$ for all $Z,W\in \H$. In order to express these
definitions in components, fix a frame $Z_1,...,Z_n$ of
holomorphic  tangent vector fields. Let $h_{\a\bb}=h(Z_\a,
Z_{\bb})$ and $g_{\a\bb}= g(Z_\a, Z_{\bb})$. The Levi curvature of
$M$ is
\begin{equation}
 \label{HA}
   H  = \frac 1 n h_\a^\a.
\end{equation}
The surface $M$ is Levi umbilical if $h_{\a\bar\b}= H
g_{\a\bar\b}$. Observe that the relation  between the Levi
curvature $H=H_\C$ and the standard mean curvature $H_\R$ is $
(2n+1)H_\R = 2n H_\C + h(T,T)$.

It is useful to compute the Levi curvature by means of a defining
function. Let $M=\{z\in \C^{n+1}:F(z)=0\}$ for some smooth
function $F:\C^{n+1}\to\R$. The holomorphic unit normal is
\begin{equation}
 \label{hollo}
    N=\sqrt 2 \frac{F_{\bar h} }{|\p F|} \p_h,\qquad\textrm{where}
    \qquad
      |\p F|^2= F_h F_\hh.
\end{equation}
The complex Hessian  $D^2 F$ induces a hermitian form on
holomorphic vector fields of $\C^{n+1}$ by letting $D^2F(U,\bar
V)= U^h V^{\bar k} F_{h\bar k}$, where  $U=U^h \p_h$ and
$\VV=V^{\kk}\p_\kk$. As observed in
 \cite{ML}, the Levi form can be written as
\begin{equation}
\begin{split}
 \label{Hesse}
   h(U,\bar V) = \frac{1}{2|\p F|} D^2F (U,\bar V),\quad U,V \in
   \H.
\end{split}
\end{equation}
Moreover, the Levi curvature of $M$ is
\begin{equation}
  \label{LeviF}
\begin{split}
 H= \frac{1}{n |\p F|} \left(
        F_{h\bar h} -
             \frac{F_{k} F_{\bar h} F_{h\bar k } }
        {|\p F|^2}  \right).
\end{split}
\end{equation}

We briefly check \eqref{Hesse}. By \eqref{NT} and \eqref{hollo},
we have
\[
    h(U,\bar V)
         = g( U, D_\VV\nu)
         = \frac{1}{\sqrt 2} g\big( U, D_\VV(N+\bar N)\big)
        = g \left( U, D_\VV
          \left(\frac{F_{h} }{|\p F|}\p_\hh\right)
          \right).
\]
As $g(F_h \p_\hh,U)=0$,  we get
\[
  g \left( U, D_\VV
          \left(\frac{F_{h} }{|\p F|}\p_\hh\right)
          \right) =
   \frac{1}{|\p F|} g( U,
         D_\VV \left(F_{h} \p_\hh\right) )
     = \frac{1}{2|\p F|}   D^2 F (U,\bar V).
\]

In order to prove \eqref{LeviF}, assume, for instance,
$F_{n+1}\neq 0$ near a point $P\in M$ and consider the local
holomorphic frame near $P$
\begin{equation}
 \label{zeta}
         Z_\a = \p_{\a}-\frac{F_\a}{F_{n+1}} \p_{n+1},\quad
         \a=1,...,n.
\end{equation}
The application of \eqref{Hesse} to the $Z_\a$'s gives
\[
  h_{\a\bb} =\frac{1}{2|\p F|}
  \Big\{F_{\a\barb } -
\frac{F_{\barb}}{F_{\overline{n+1}}} F_{\a\overline{n+1}}
-\frac{F_\a}{F_{n+1}}F_{n+1,\barb} + \frac{F_\a
F_\barb}{|F_{n+1}|^2} F_{n+1,\overline{n+1}}\Big\}.
\]
The metric tensor and its inverse are respectively
\[
       g_{\a\bar\b}  = \frac 12 \left(
                    \d_{\a\b}+\frac{F_\a F_{\bar\b} }{
                    |F_{n+1}|^2}\right)
                        \quad\text{and}\quad
       g^{\a\bar\b} =2
                 \left(
                    \d_{\a\b}-\frac{F_\b F_{\bar\a} }{
                    |\p F|^2}\right).
\]
Then, a short computation gives
\[
\begin{split}
 H  = \frac{1}{n}  g^{\a\bar\b}h_{\a\bb}=
           \frac{1}{n |\p F|} \left(
        F_{h\bar h} -
             \frac{F_{k} F_{\bar h} F_{h\bar k } }
        {|\p F|^2}  \right).
\end{split}
\]

In the next proposition we collect some useful identities.

\begin{proposition}
 \label{popolo}
Let $M\subset\C^{n+1}$ be an oriented surface with real unit
normal $\nu$, $T=J(\nu)$ and holomorphic unit normal $N$. Then:
\begin{itemize}
  \item[i)] $g(D_Z N,\NN) = g([T,Z],T)$ for all $Z\in \Gamma(\H)$;
  \item[ii)] $g(D_Z N,\NN) = i h(T,Z)$ for all $Z\in\Gamma(\C T
  M)$;
  \item[iii)] $g([Z,\WW],T)=- 2 i h(Z,\WW)$ for all
  $Z,W\in\Gamma(\H)$.
\end{itemize}
\end{proposition}

\begin{proof}
Note that $g(D_Z T,  T) = 0$, because $T$ is real. Moreover, by
\eqref{NT}, we have for any $Z\in\Gamma(\H)$
\[
\begin{split}
   g([T,Z], T) & = g( D_T Z- D_Z T,T)
     %\\&
   = \frac 12 g(  D_{N-\bar N} Z,\bar N \rangle.
 % \\&
 % =  \frac 12 g(D_N Z, \NN)
 %   -\frac 12 g(D_{\bar N} Z, \NN).
\end{split}
\]
We used the orthogonality  $g(D_{N-\bar N} Z,\NN)=0$, which holds
because $Z$ is holomorphic. Thus
\[
\begin{split}
     2 g( [T,Z], T)  &
      = g( D_Z N +[N,Z], \NN)
       -  g( D_Z\bar N +[\bar N, Z],\NN)
   \\&
   =  g( [N-\bar N, Z],\NN)+g(D_Z N,\NN)
   \\&
   = \frac{1}{ i}  g([T,Z] , \nu+ i T)+g( D_Z N,\NN)
   =  g([T,Z],T) + g(D_Z N,\NN).
\end{split}
\]
We used again  \eqref{NT} and $g([T,Z],\nu)=0$. This proves i).

In order to check ii), note that
\[
\begin{split}
     g(D_Z N,\NN) & =g(D_Z N, N+\bar N)
                    = \sqrt 2 g( D_Z N, \nu)
                    =  g( D_Z \nu, \nu)
                     -i g(D_Z T, \nu)
                \\& = -i g(D_Z T, \nu)
                    = i h(T, Z).
\end{split}
\]

Identity iii) is proved in \eqref{campanella}--\eqref{campanello}.
\end{proof}

Now we discuss a couple of examples showing  the existence of non
trivial Levi umbilical surfaces.

\begin{example}[Boundary of spherical tubes]
 \label{esso}
The surface $ M =\{ z\in \C^{n+1} :F(z)=0\}$, where
\[
     F(z)=\frac 12 \sum_{h=1}^{n+1} (z_h+\bar z_h)^2-1,
\]
is a Levi umbilical cylinder with spherical section having
constant Levi curvature. Indeed, the complex derivatives of $F$
are $F_h=F_{\bar h} = z_h+\bar z_h$ and $F_{h\bar k}= \d_{hk}$.
Then $|\p F|=\sqrt 2$ and, by \eqref{LeviF}, the Levi curvature is
$H=1/\sqrt 2$. The complex Hessian of $F$ is the identity and, by
\eqref{Hesse}, the condition $h_{\a\bb}= \frac{1}{\sqrt 2}
g_{\a\bb}$ is identically satisfied on $M$.
\end{example}

\begin{example}
 \label{es2}
 It is possible to construct compact Levi umbilical surfaces by
 polyharmonic perturbations of the sphere. Consider
\begin{equation}
 \label{Emme}
         M=\left\{z\in\C^{n+1} : |z|^2+\lambda\Phi(z) =1\right\},
\end{equation}
where $\la$ is a real parameter and
\[
      \Phi(z) = \frac 12 \sum_{h=1}^ {n+1}(z_h^2+\bar z_h^2).
\]
The derivative of the defining function  $F(z)=|z|^2+\la \Phi
(z)-1$ are  $F_h=\bar z_h+\la z_h$ and $F_{h\bar k}=\d_{hk}$. On
the set $M$ we have
$
    |\p F(z)|^2 = 2-(1-\la^2)|z|^2.
$
Then,  $|\p F|$ is constant on $M$ if and only if $\la=0,1,-1$. If
$|\la|<1$, $M$ is a smooth compact surface bounding the region
$\{z\in\C^{n+1}: F(z) <0 \}$. Indeed, $M$ is an ellipsoid: letting
$z=x+iy$, we have $F(z)=(1+\la)|x|^2 +(1-\la)|y|^2-1$. Moreover,
on $M$
$
     |\p F(z)|^2
       = 2-(1+\la) (1-2\la |x|^2)\geq 1-\la>0.
$
By formula \eqref{LeviF}, the Levi curvature of $M$ is
$
       H={|\p F|^{-1}}.
$
The complex Hessian of $F$ is the identity and, by  \eqref{Hesse},
the surface $M$ is Levi umbilical and
\begin{equation}
\label{acca}   h_{\a\bb}= \frac{1}{|\p F|} g_{\a\bb}.
\end{equation}

Many other examples of compact Levi umbilical surfaces can be
constructed, taking as $\Phi$ in \eqref{Emme} any polyharmonic
function, i.e.~any smooth function satisfying $\Phi_{h\bar k}=0$.
In fact, the complex Hessian of the corresponding defining
function is the identity. Therefore  condition \eqref{acca} is
satisfied.

\end{example}

\section{The connection and its properties}
\label{cnx} \setcounter{equation}{0}

In this section, we define the covariant derivative  $\D$ on an
oriented, smooth hypersurface $M\subset\C^{n+1}$ starting from the
standard connection $D$ in $\C^{n+1}$. A vector field $V\in
\Gamma( \C T M)$ can be uniquely decomposed as
\begin{equation}
  \label{deco}
      V = Z +\bar W + f T,
\end{equation}
where $Z, W \in \Gamma(\H )$ and $f\in C^{\infty}(M)$ is a complex
valued function.  We define $\nabla :\Gamma(\C TM)\times \Gamma
(\C TM)\to\Gamma(\C TM)$ by letting, for $U,V\in \Gamma(\C TM)$
with $V$ as in \eqref{deco},
\begin{equation}
\label{lallp}
   \D_U V = D_U Z - g(D_U Z, \NN) N
          + D_U \bar W - g( D_U \bar W, N)\bar N
          + (U f) T.
\end{equation}
Here, $N$ is the holomorphic unit normal. Equivalently, let for
$U\in \Gamma(\C TM)$ and $Z,W\in \Gamma(\H)$
\begin{equation}
\begin{split}
\label{lillo}
   &
   \D_U Z = D_U Z - g(D_U Z, \NN) N,\\&
   \D_U\bar W = D_U \bar W -g( D_U \bar W,
         N) \bar N, % \quad \text { and}
   \\&
   \D_U T=0 .
\end{split}
\end{equation}

We have the following
\begin{theorem} \label{gurzo}
 $\D$  is a complex  connection  on $M$
 and satisfies the following properties:
\begin{itemize}

\item[(C1)]
  $\overline{\D_U V}= \D_{\UU}\VV$ for all $U,V\in\Gamma (\C TM)$;

\item[(C2)]
  $\D_U (J(V))=J(\D_U V)$ for all  $U, V\in \Gamma( \C TM)$;

\item[(C3)] The bundles $\H$ and $\HH$ are parallel;

\item[(C4)] $\D g=0$;

\item[(C5)] $\mathrm{Tor}_\D(U,V)=0$ for all $U, V\in \Gamma(\H)$;

\item[(C6)]
$  \mathrm{Tor}_\D(U,\VV)=-g( [U,\VV],T) T$ for all $U,
V\in\Gamma(\H)$.

\end{itemize}
\end{theorem}

\begin{proof}
Properties (C1), (C2) and the fact that $\D$ is a connection are
easy and we omit their proof.

Property (C3) amounts to say that the covariant derivative of a
holomorphic (resp.~antiholomorphic) vector field is still a
holomorphic (resp.~antiholomorphic) vector field. But this   is an
immediate consequence of definition \eqref{lallp} and of the
orthogonal decomposition
$
    T^{1,0}_P \C^{n+1} = \H_P \oplus \C N_P,
$
at any point $P\in M$.

In order to prove property (C4),  let
\[
 V_1= Z_1 + \WW_1 + f_1 T,\qquad
 V_2= Z_2 + \WW_2 + f_2 T,
\]
where $Z_1,Z_2, W_1, W_2\in \Gamma(\H)$ and $f_1, f_2$ are complex
valued functions. By the metric property of the standard
connection $D$ in $\C^{n+1}$, we have
\begin{equation}
\begin{split}
\label{spqr}
   U g( V_1, \VV_2) &
      = U g(Z_1, \ZZ_2)
      + U g(\WW_1, W_2)
      + U g(f_1 T, \bar f_2 T)
\\&
      = g( D_U Z_1, \ZZ_2)
      + g(Z_1, D_U \ZZ_2)
      + g (D_U\WW_1, W_2)
      + g (\WW_1, D_U W_2)
\\& \qquad
      + g(D_U(f_1 T), \bar f_2 T)
      + g (f_1 T, D_U (\bar f_2 T)).
\end{split}
\end{equation}
We claim that the following identities hold
\begin{equation}
\label{idill}
\begin{array}{ll}
 g( D_U Z_1, \ZZ_2) =  g(\D_U Z_1, \VV_2),      &
 g(D_U\WW_1, W_2) = g ( \D_U \WW_1, \VV_2),
 \\
 g( Z_1, D_U \ZZ_2) = g( V_1, \D_U \ZZ_2),  &
 g(\WW_1,D_U W_2)= g(V_1 ,\D_UW_2).
\end{array}
\end{equation}
We check the first one only. Since $g( N, \ZZ_2)=0$,  we have
$g(D_U Z_1, \ZZ_2)= g(\D_U Z_1, \ZZ_2)$  and  property (C3) gives
$g(\D_U Z_1, \ZZ_2) = g(\D_U Z_1, \VV_2)$. The following
identities also hold
\begin{equation}
 \label{dillo}
  g(D_U(f_1 T), \bar f_2 T) =g( (Uf_1) T ,\VV_2),
     \qquad
   g( f_1 T, D_{U}(\bar f_2T)) = g( V_1,(U \bar f_2) T).
\end{equation}
We check the first one. Since $g( D_U T , T)=0$, then  $g( D_U(f_1
T), \bar f_2 T) =g( (Uf_1) T , \bar f_2 T) + g(f_1 D_U T , \bar
f_2 T) = g( (Uf_1) T , \bar f_2 T)$.   But  $T$ is orthogonal to
$Z_2$ and $\bar W_2$. Thus we get the claim. Replacing
\eqref{idill} and \eqref{dillo} into \eqref{spqr} we get $U g(V_1,
\VV_2) = g(\D_U V_1, \VV_2) + g(V_1, \D_U \VV_2)$, which means $\D
g=0$.

Statement (C5), $\mathrm{Tor}_\D (U,V)=0$ for $U,V\in \Gamma(\H)$,
follows from $\mathrm{Tor}_D (U,V)=0$ and $[U,V]\in\Gamma(\H)$.
Concerning property (C6),  observe that a connection leaving $\H$
and $\HH$ parallel cannot be, in general, torsion free, because
the horizontal distribution needs not be integrable (in other
words, it may be $[\H,\HH]\nsubseteq \H\oplus \HH$). Take  $W\in
\Gamma(\mathcal H\oplus\HH)$. Then
\[
   g( \mathrm{Tor}_\D(U,\VV), \WW) =
   g\big (\mathrm{Tor}_D(U,\VV)
      +g(D_\VV U,\NN) N
      -g(D_U\VV,N) \bar N,\WW\big) =0,
\]
because $ \mathrm{Tor}_D(U,\VV)=0$ and $g(N,\WW) = g(\bar
N,\WW)=0$. Then $ \mathrm{Tor}_\D(U,\VV)=\lambda T$ for some
function $\lambda$ and, by (C3),
$
   \lambda= g( \mathrm{Tor}_\D(U,\VV),T) =
   -g( [U,\bar V],T).
$
\end{proof}

Recall that the restriction of the hermitian product in $\C^{n+1}$
to $\C TM$ induces the orthogonal decomposition
\begin{equation}
 \label{dc}
      \C TM= \H\oplus \HH \oplus \C T.
\end{equation}
Denote by $\Pi_\H:  \C TM\to \H$ the projection onto $\H$ and by
$\Pi_\HH$  the projection onto $\HH$. Then it is easy to check
that for $U,V\in \Gamma(\H)$, we have
\begin{equation}
 \label{comma}
   \D_\UU V = \Pi_\H([\UU,V])
   \quad
   \textrm{and}
   \quad
   \D_U \VV = \Pi_\HH([U,\VV]).
\end{equation}
This follows from (C3) and (C6).

\begin{remark}\label{webb} If $M$ is a strictly pseudoconvex CR manifold,
there is a natural connection associated with a given contact form
$\theta$, which was introduced by Tanaka and Webster in \cite{T}
and \cite{W}. Although it was designed for different scopes from
ours, we highlight some analogies and differences between our
connection $\D$ and the Tanaka-Webster one.

The Levi form $(Z,W)\mapsto -i d\theta(Z,\WW)$ is a non degenerate
hermitian form on $\H$.   Then, $\theta$ induces a decomposition
of $\C TM$ similar to \eqref{dc}. The vector field $T$ is replaced
in this construction by the {\it characteristic vector field}
$T'$, defined by $\theta(T')=1$ and $d\theta(T, Z)=0$ for all
$Z\in \H$. In the Tanaka-Webster connection, the Levi form
$d\theta$ essentially plays the role of the metric and is required
to be parallel. Covariant derivatives of holomorphic vector fields
along  antiholomorphic ones are defined by relations analogous to
\eqref{comma} (see \cite[Lemma 3.2, p. 31]{T}), but with the
projections $\Pi'_\H$ and $\Pi'_{\bar \H}$ induced by $d\theta$.
The characteristic vector field $T'$ of $\theta$ is in general
different from $T$ for any choice of the contact form $\theta$.

 Similarly to the Tanaka-Webster connection, the property $\D
T=0$ is  forced by (C4) and (C5). Indeed, the one dimensional
bundle generated by $T$ is the orthogonal complement with respect
to the parallel metric $g$ of the parallel bundle $\H\oplus\H$.
Then $\D_U T =\lambda T$ for some function $\lambda$ and $U\in \C
TM$. But, since $T$ is real,
$
  0=U g(T, T)
   = 2g( \D_U T, T=2\lambda
$.
 Therefore $\D T=0$.
\end{remark}

\begin{remark}
  The connection $\nabla$ is not uniquely determined on the whole
tangent bundle $\G(\C TM)$ by properties (C1)--(C6). In
particular, $\D_T U$ with $U\in\G(\H)$ is not uniquely determined.
In \eqref{lillo}, we let  $\D_T U= D_T U-g( D_T U, \NN)N$. An
alternative possibility, consistent with \eqref{comma}, is to set
\[
   \D'_T U = \Pi_\H ([T,U])
   \quad
   \textrm{and}
   \quad
   \D'_T \UU = \Pi_{\bar \H} ([T,\UU]).
\]
The resulting connection $\D'$ still satisfies (C1)--(C6). Our
choice  $\D$, however, seems to be more suitable than $\D'$ to
work with Codazzi equations.

\end{remark}

\begin{remark}
  \label{snam}
The real tangent bundle has the orthogonal decomposition $TM=
\mathrm{Re}(\H\oplus\HH) \oplus \R T$. Then, for $Y\in
\Gamma(\mathrm{Re}(\H\oplus\HH))$, $V\in \Gamma(TM)$ and $f$ real
function, we have
\begin{equation}
 \label{gillo}
  \D_V (Y+f T)= D_V Y - g( D_V Y, \nu)\nu
       -g( D_V Y, T) T + (Vf) T.
\end{equation}
%\end{equation}
Indeed, taking $X=Z+\bar Z$ with holomorphic $Z$, we have
%\begin{equation}
 % \label{uuu}
\[
 \begin{split}
   \D_V(Z+\ZZ) & =\D_V Z + \D_V \ZZ =  D_V Z - g( D_V Z, \NN)N
    +  D_V \ZZ - g( D_V \ZZ, N)\NN
\\&
 = D_V X - g(D_V Z, \nu) (\nu-iT) - g( D_V \ZZ, \nu  ) (\nu+i T)
\\&
 =D_V X - g( D_V X, \nu)\nu
  +i g(D_V(Z-\ZZ), \nu) T
\\& =D_V X - g(D_V X, \nu)\nu + g( D_V(J(X)), -J(T))T
%\\&
% =D_V X - \langle D_V Z, \nu
%\rangle\nu +i \langle D_V(Z-\ZZ), \nu\rangle T
\\& =D_V X - g (D_V X, \nu)\nu -g(D_V X, T)T.
\end{split}
\]
%\end{equation}
We used $\nu=-J(T)$ and the property $D_V \circ J = J\circ D_V $.
\end{remark}

\section{Codazzi equations}
\setcounter{equation}{0}

In this section we compute the system of Codazzi equations.

\begin{theorem}
  \label{coda}
The Levi form $h$ on a hypersurface $M\subset\C^{n+1}$ satisfies
the following Codazzi equations \begin{subequations}
\begin{eqnarray} \label{azzo}
 &\nabla_\b h_{\a\barga}- \nabla_\barga h_{\a\b}
  =  i h_{\a\b} h_{\gg 0}
    -i h_{\a\barga} h_{\b 0}
    -2i h_{\b\gg} h_{\a0},
\\
  \label{gotta}
  &\D_\bb h_{\a 0}- \D_0 h _{\a\bb} =
        i h_{\a\la} h^\la_\bb
      - i h_{\a\ll} h^\ll_\bb
      + i h_{\a\bb} h_{00},
\\
 \label{geog}
  &\D _\b h_{\a0} - \D_0 h_{\a\b} =
         i h_{\b\ll}h^\ll_\a
      -  i h_{\b\la} h^\la_{\a}
      + i h_{\a\b} h_{00}
      - 2i h_{\a0} h_{\b0},
\\
     \label{noni}
  &  \D_\a h_{00} - \D_0 h_{\a 0} =
        2 i h_\a^\la h_{\la0}
       -2 i h_\a^\ll h_{\ll 0}
      - i h_{\a\la}h^\la_0
      + i h_{\a\ll} h^\ll_{0}
      - i h_{\a0} h_{00}.
\end{eqnarray}
\end{subequations}
\end{theorem}

\begin{proof}
The proof relies on the fact that the standard connection $D$ in
$\C^{n+1}$ has vanishing curvature. We shall also use several
times the formula
\begin{equation}
\label{hh}
       D_Z U= \D_Z U -\sqrt 2 h(U,Z) N,\qquad
       U\in \G(\H),\,
       Z\in \G(\C T M).
\end{equation}
Let $Z,W\in \G(\C TM)$ and $ U\in \G(\H)$. Denote by $R_D$ and
$R_\D$ the standard curvature endomorphisms of $D$ and $\D$. Using
\eqref{hh}, we have
\begin{equation}
     \label{pistola}
\begin{split}
   0 & =  R_D(Z,W)U  =  D_Z D_{W} U - D_{W} D_ZU
          - D_{[Z,W]} U
\\&
       =  D_Z\big( \D_{W}U - \sqrt 2 h(U, W) N \big)
       -  D_{W} \big( \D_{Z}U - \sqrt 2 h(U, Z) N \big)
% \\& \qquad
      - D_{[Z,W]} U
     % -  \D_{[Z, \bar W]}U +\sqrt 2 h(U, [Z,\WW]) N
\\&
       =  R_\D(Z,W) U
       - \sqrt 2 h(U, W) D_Z N +  \sqrt 2 h(U,  Z) D_{W}N
\\&\quad
       - \sqrt 2 \Big(  Z h(U, W) - h(\D_Z U, W)
       - W h(U, Z) + h(\D_{W} U, Z)
       - h(U, [Z,W]) \Big) N.
\end{split}
\end{equation}
Multiplying by $\NN$ and using $g(  R_\D(Z,W) U,\NN) =0$, we get
the equation for $h$
\begin{equation}
  \label{pisso}
  \begin{split}
   \D_ Z h (U,W)
  -\D_W h(U,Z)=
   h(U, \mathrm{Tor}_\D (W,Z) )       &
  - h(U, W) g( D_Z N,\NN) \\ &
  + h(U, Z) g(D_W N,\NN) ,
\end{split}
\end{equation}
where $
     \nabla_Z h(U, W)
     = Z h (U, W) - h(\nabla_Z U,W)-h(U,\nabla_Z W)
     $
is the covariant derivative of $h$. Note that by Proposition
\ref{popolo}, we have $g(  D_Z N,\NN) = i h(Z,T)$ for any $Z\in
\G(\C TM)$.

In order to prove \eqref{azzo},  take $Z,W,U\in \G(\H)$ and write
\eqref{pisso} with $\WW$ instead of $W$. By Theorem \ref{gurzo},
the torsion satisfies $\mathrm{Tor}_\nabla(Z,\bar W)=-g( [Z,\bar
W],T) T$. Moreover, by Proposition \ref{popolo} we have $g(
[Z,\bar W],T) = - 2i h(Z,\bar W)$.  Thus, equation \eqref{pisso}
becomes
\begin{equation}
  \label{pisso1}
  \nabla_Z h (U,\bar W)- \nabla _{\bar W}h(U,Z) =
     i h(U,Z) h(\WW,T)
  -  i h(U,\bar W) h(Z,T)
   -   2i h(Z,\bar W) h(U,T).
\end{equation}
This is formula \eqref{azzo}.

In order to prove \eqref{gotta}, we take $Z,U\in \G(\H)$. By
\eqref{pisso}, we have
\begin{equation}
\begin{split}
 \label{gui}
 \D_{\bar Z} h (U,T)-\D_T h(U,{\bar Z})= h(U,\D_T{\bar Z}+[{\bar Z},T])
             &      - i h(U,T) h(\bar Z,T)
             \\ &
                    + i h(U,{\bar Z}) h(T,T),
\end{split}
\end{equation}
because $\D T=0$.

We analyze the right hand side of \eqref{gui}. By
$\mathrm{Tor}_D({\bar Z},T)=0$ and the second equation of
\eqref{lillo}, we have
\begin{equation}
 \label{vino}
\begin{split}
 \D_T{\bar Z}+[{\bar Z}, T]&= D_{\bar Z} T
             - g( D _ T{\bar Z}  , N)\NN
   %     \\ &
   %               = D_{\bar Z} T - \langle D_T {\bar Z},N+\NN\rangle \NN
        \\ &
                  = D_{\bar Z} T - g( D_{\bar Z} T,\nu) \nu
                            -i g(  D_T {\bar Z},\nu)T.
\end{split}
\end{equation}
We also used $g(  [{\bar Z},T],\nu) =0$, which implies $g( D_\ZZ
T, \nu)=g(  D_T \ZZ , \nu)$.  The vector field $V=D_{\bar Z} T -
g( D_{\bar Z} T,\nu) \nu$ is tangent to $M$ and $g( V,T) =0$.
Therefore, for any holomorphic frame $Z_1,...,Z_n$, we have
\begin{equation}
 \label{olio}
\begin{split}
     V & = g^{\lambda\mm} g(  D_{\bar Z} T, Z_{\mm})Z_\lambda
        +g^{\mu\ll} g( D_{\bar Z} T,Z_\mu )Z_\ll \\
       & = i g^{\lambda\mm} g(  D_{\bar Z} \nu, Z_\mm) Z_\lambda
        - i g^{\mu\ll} g(  D_{\bar Z} \nu,Z_\mu) Z_\ll \\
       & = i g^{\lambda\mm}  h(Z_\mm, {\bar Z}) Z_\lambda
        - i g^{\mu\ll} h(Z_\mu,{\bar Z}) Z_\ll   .
\end{split}
\end{equation}
In order to get the second equality in \eqref{olio}, we used the
isometry $J$ and the relations $T=J(\nu)$, $J(D_{\bar Z}
T)=D_{\bar Z}( J(T))$, $J(Z_\mu)=i Z_\mu$ and $J(Z_\mm)=-i Z_\mm$.
Thus,  \eqref{vino}--\eqref{olio} give
\begin{equation}
 \label{bios}
\begin{split}
     h(U,\D_T{\bar Z}+[{\bar Z},T])
     = i h(T,{\bar Z}) h(U,T)
     & + i g^{\lambda\mm} h(Z_\mm,{\bar Z}) h(U,Z_\lambda)
   \\&
      -i g^{\mu\ll} h(Z_\mu,{\bar Z}) h(U,Z_\ll).
\end{split}
\end{equation}
Replacing \eqref{bios} into \eqref{gui}, we finally find
\[
\begin{split}
   \D_{\bar Z} h(U,T) -\D_T h(U,{\bar Z}) = i h(U,{\bar Z}) h(T,T)
       & + i g^{\lambda\mm} h(Z_\mm,{\bar Z}) h(U,Z_\lambda)
       \\&
      -i g^{\mu\ll} h(Z_\mu,{\bar Z}) h(U,Z_\ll),
\end{split}
\]
which is identity \eqref{gotta}.

In order to prove \eqref{geog},  take $Z,U\in \G(\H)$. By
\eqref{pisso}, we have
\begin{equation}
\begin{split}
 \label{venusia}
    \D_Z h(U,T) -\D_T h(U,Z) & = i h(U,Z)h(T,T)-i h(U,T)h(Z,T)
  \\&\qquad                          +h(U,\D_T Z + [Z,T]).
\end{split}
\end{equation} On conjugating \eqref{vino}, we find $\D_T Z
+[Z,T] = \bar V - i h(Z,T) T$, where the vector field $\bar V= D_Z
T -g(  D_Z T,\nu)\nu$ is, by \eqref{olio},
\[
    \bar  V= i g^{\la\mm} h (Z, Z_\mm) Z_\la
       -i g ^{\mu\ll} h ( Z, Z_\mu) Z_\ll.
\]
Thus, equation \eqref{venusia} reads
\[
\begin{split}
     \D_Z h(U,T)-   \D_T h(U,Z) = & i h(U,Z)h(T,T)
                  -2 i h(U,T)h(Z,T)
             \\
             & +i g^{\la\mm} h(U,Z_\la) h (Z, Z_\mm)
             -i g ^{\mu\ll} h ( Z, Z_\mu) h(U, Z_\ll).
\end{split}
\]
The proof of identity \eqref{geog} is accomplished.

In order to prove \eqref{noni},  take $Z\in \G(\H)$ and  start
from the identity
\begin{equation}
 \label{pio}
       D_T D_Z T - D_Z D_T T - D_{[T,Z]} T=0.
\end{equation}
Observe that  $D_Z T = U  -h(Z, T) \nu$ for some $U\in
\G(\H\oplus\bar\H)$, because $g(  D_Z T, T) =0$. Precisely, as in
\eqref{olio}, we have
\begin{equation}
 \label{uva}
     U = i g^{\la\mm} h(Z, Z_\mm) Z_\la - i g^{\mu\ll} h(Z, Z_\mu)
     Z_\ll.
\end{equation}
Then
$
        D_T D_Z T = D_T U  - Th(Z,T) \nu - h(Z,T) D_T \nu,
$
and multiplying by $\nu$,
\begin{equation}
 \label{pia}
     g(     D_T D_Z T, \nu)
      =- h(U,T) -Th(Z,T),
\end{equation}
because $g(  D_T\nu,\nu)=0$.

We analyze the second term in the left hand side of \eqref{pio}. A
computation similar to \eqref{olio} furnishes
\begin{equation}
 \label{Vio}
    D_T T-g(  D_T T,\nu)\nu =
      i g^{\la\mm} h(Z_{\mm},T) Z_\la
     -i g^{\mu \ll} h(Z_\mu,T) Z_\ll = W,
\end{equation}
where $W\in \G(\H\oplus \bar\H)$ is defined by the last equality.
Thus,
\begin{equation}
 \label{pii}
   g(   D_Z D_T T,\nu)
      =  -h(Z, W)  - Z h(T,T).
\end{equation}

Finally, we study the third term in the left hand side of
\eqref{pio}. We have
\[
\begin{split}
 [T,Z] & = D_T Z- D_Z T
         = \D_T Z + g(  D_T Z,\NN) N - D_Z T \\&
         = \D_T Z -h(Z,T) (\nu-i T) - D_Z T
         = \D_T Z - U + i h(Z,T) T,
\end{split}
\]
where $U$ is defined after \eqref{pio}. This yields
\begin{equation}
 \label{piu}
  g(  D_{[T,Z]} T,\nu) = - h([T,Z],T)
%                    = - h(\D_T Z - U + i h(Z,T) T,T)
                    = - h(\D_T Z,T) + h(U ,T)
                      -  i h(Z,T)h(T,T).
\end{equation}
Multiplying \eqref{pio} by $\nu$ and using \eqref{pia},
\eqref{pii} and \eqref{piu},
 we obtain
\[
  \D_Z h(T,T)-\D_ Th(Z,T) =  2h(U ,T) - h(Z, W) - i
  h(Z,T)h(T,T).
\]
Replacing the expressions for $U$ and $W$ in \eqref{uva} and
\eqref{Vio}, we get formula \eqref{noni}.
\end{proof}

\begin{remark}
 \label{remi}
 The second fundamental form $h$ satisfies also other Codazzi
equations. For instance, we have
\begin{subequations}
\begin{eqnarray}
 &
   \nabla_\a h_{\b\barga} -  \nabla_\b h_{\a\barga}   =
   i   h _{\b\barga}h_{\a0}-  i h _{\a\barga} h_{\b0},
   \label{azzo2}
\\&
\label{azzo3}
    \nabla_\a h_{\b\ga}    - \nabla_\ga h_{\b\a}=
    i  h_{\b\a} h_{\ga0} - i h_{\b\ga} h_{\a0},
%   \\&
%\label{azzo4}
%   \nabla_{\bar\a} h_{\b\barga} - \nabla_{\barga} h_{\b\bar\a}=
%    i h_{\b\bar\a} h_{\barga0}-i  h_{\b\barga} h_{\aa0}.
\end{eqnarray}
\end{subequations}
Identity \eqref{azzo2} can be obtained interchanging $\a$ and $\b$
in identity \eqref{azzo} and taking the difference of the two
equations. Identity \eqref{azzo3} follows from \eqref{pisso} on
choosing $Z,U,W\in \G( \H )$ and using $\mathrm{Tor}_\D(W,Z)=0$.

Notice also that, letting $Z,U,V,W\in\G(\H)$ and multiplying
identity \eqref{pistola} by $V$, we get the Gauss-type equation
\[
   g( R_{\D}(Z, \bar W) U, \VV)
   = 2\big\{ h(U, \bar W)h(\bar V, Z) - h(U,Z) h (\bar V, \bar
   W)\big\}.
\]
\end{remark}

\section{Classification results}
\setcounter{equation}{0}

In this section we prove the following results:
\begin{theorem}
       \label{teor}
Let $M\subset \C^{n+1}$, $n\ge 2$, be a $(2n+1)$-dimensional,
connected  Levi umbilical surface with constant Levi curvature
$H\neq 0$.  Then $M$ is contained either in a sphere or in the
boundary of a spherical tube.
\end{theorem}

\begin{theorem}
\label{klino} Let $M$ be a connected pseudovonvex hypersurface in
$\C^{n+1}$, $n\geq 1$, with constant Levi curvature $H \neq 0 $
and $h_{\a\b}=0$. Then, up to a complex isometry, $M$ is a
contained in a sphere  or in a cylinder of the form
\begin{equation}
 \label{cyllo}
       \Big\{ z\in\C^{n+1} : \sum_{i=m}^{n+1} |z_i|^2
       =r^2\Big\},\quad r>0, \quad 1\leq m\leq n.
\end{equation}
\end{theorem}

\begin{remark}
  \label{ria}
The only compact surface among the ones defined in \eqref{cyllo}
is the sphere.  Theorem \ref{klino}  improves \cite[Theorem
5.2]{K}, because  we assume neither compactness nor strict
pseudoconvexity of $M$.

A slight modification of the argument also shows that if strict
pseudoconvexity (but not compactness) is added as hypothesis in
Theorem \ref{klino}, then the surface $M$ must be contained in a
sphere.
\end{remark}

\begin{proof}[Proof of Theorem \ref{teor}]
Possibly changing the orientation of $M$,  assume $H>0$.
  Observe preliminarily that, given an
orthonormal frame $Z_\a$, by Proposition \ref{popolo} part iii),
we have
\begin{equation}
 \label{rnk}
    \sum_{\a=1}^n g([Z_\a, Z_\aa], T)= -2i n H \neq 0,
\end{equation}
provided that $H\neq 0$. Then at least one term in the sum is non
zero and the distribution $\mathrm{Re}(\H\oplus\HH)$ is bracket
generating.

We accomplish the proof in several steps.

\smallskip \noindent\it Step 1. \rm    We claim
 that
 \begin{equation}
 \label{carro}
 h_{\a0}=0.
 \end{equation}
Indeed, contracting the indices $\a$ and $\gg$ in the Codazzi
equation \eqref{azzo2}, we get
\begin{equation}
 \label{stachel}
   \nabla_\a  h_{\b}^\a - \nabla_\b h_\a^\a=
    i  h _\b ^\a  h_{\a0} -  i  h _\a^\a h_{\b0}.
\end{equation}
The fundamental form satisfies $h_{\a\bar\b}= H g_{\a\bar\b}$ and
thus $h_\a^\b = H \d_\a^\b$ and $h_\a^\a=n H$.   Then the left
hand side in \eqref{stachel} vanishes.
 Therefore
$       (n-1) H  h_{\b0}=0. $ The claim follows.

As a consequence of \eqref{carro}, it turns out that $h$ satisfies
the identities
\begin{subequations}
\begin{eqnarray}
\label{gotto} &         h_{\a}^\bb h_\bb^\mu =(H^2
-h_{00}H)\delta_\a^\mu,
  \\
\label{geografo}
  &         \D_0 h_{\a\b}+ih_{00}h_{\a\b}=0, \\
\label{nonino}
  &       \D_\a h_{00}=0.
\end{eqnarray}
\end{subequations}
To show \eqref{gotto}, observe that, since  \eqref{carro} holds
and $h_{\a\bb}=
 Hg_{\a\bb}$,
the left hand side of identity \eqref{gotta} vanishes. Thus, using
again
 $h_{\a\bb}=H g_{\a\bb}$ in the right-hand side, we find the
equation
\[
    h_{\a\la} h^\la_\bb =
        (H^2 - h_{00} H) g_{\a\bb}.
\]
Contracting with $g^{\mu\bb}$ yields \eqref{gotto}. Equations
\eqref{geografo} and \eqref{nonino} follow from \eqref{geog} and
\eqref{noni}, letting $h_{\a0}=0$ and $h_{\a\bb}=H g_{\a\bb}$.

Notice  also that   equation \eqref{gotto} gives
\begin{equation}
   \label{geografia}
         |h_{\a\b}|^2:= h_\a^\bb h_\bb^\a  =
         nH(H - h_{00}),
\end{equation}
which  implies  $h_{00} \le H$.
 Moreover, equation \eqref{nonino} and $\D
T=0$ give $ Z h(T, T) = \D_Z h(T, T) =0$ on $M$ for any $Z\in \H$.
On conjugating, the equation is satisfied also for all
$Z\in\bar\H$.  Since $M$ is connected, from \eqref{rnk} it follows
that
\begin{equation}
 \label{dido}
     h(T, T)=\text{constant}=h_{00} \quad \text{on }M.
\end{equation}

Take $P\in M$ and  denote  by $L$
 the shape
operator,  $L(X)=D_X\nu$, $X\in T_PM.$

\smallskip
\noindent \it Step 2. \rm If  $X\in T_PM$ is an eigenvector of $L$
with $|X|=1$ and $g( X,T)=0$, then $Y=J(X)$ is an eigenvector of
$L$ with $|Y|=1$.

\smallskip

Indeed, assume that  $L(X)=\la X$ for some $\la\in\R$ and let
$Z=X-i J(X)=X-iY\in \H_P $. By \eqref{carro}, since $L(X)$ is
orthogonal to $T$, \[ 0=h(Z,T) = g(L(X)-i L(Y), T) = -i g(L(Y), T)
.
\]
Therefore $L(Y)$  is orthogonal to $T$.
 Moreover, by the symmetry of $L$,  $
   g( L(Y),X) = g( L(X),Y)  = \la g(
   X,Y)  = \la g(
   X,J(X))  =0$.
Finally, if  $W\in \H_P  $ satisfies $g(Z,\WW)  =0$, it must be
also $g(X, \WW) =0$ and thus $g( L(X),\WW)  = \la g( X,\WW) =0$.
Since $M$ is Levi umbilical, we also have
 $g( L(Z),\WW)  =  H g(Z,\WW) $. Eventually, we get
\[
      g(  L(Y),\WW)  = i g(  L(X)-i L(Y),\WW)
                             = i g(  L(Z),\WW)
                             = i H g(Z,\WW)=0.
\]
Taking the conjugate we also find $g( L(Y),  W)  =0$. Ultimately,
we showed that $L(Y)$ has no component orthogonal to $Y$ and our
claim is proved.

\medskip
\noindent \it Step 3. \rm At any point $P\in M$ there exists an
orthonormal basis $\{ X_\a, Y_\a = J(X_\a), T:\a=1,\dots, n\}$ of
$T_PM$ such that
\begin{equation}
\begin{split}
\label{acquaragia}
 & L(X_\a) = (H+\sqrt{H^2 - h_{00} H}) X_\a,
    \\
 & L(Y_\a) =(H-\sqrt{H^2 - h_{00} H}) Y_\a, % \quad\text{and}
    \\
 & L(T)=h_{00} T.
\end{split}
\end{equation}

\medskip

Note first that, by \eqref{carro},
    $h(T, X)=0$ for any
$X\in\H_P\oplus\bar\H_P$. Then we have $L(T) = h_{00 } T$, by
\eqref{dido}, and the orthogonal complement of $T$ at any point
$P\in M$ is an invariant subspace for $L$. We diagonalize $L$
restricted to this invariant subspace. By {\it Step 1}, for any
eigenvector $X_\a$ with eigenvalue $\la_\a$, there is an
eigenvector $Y_\a = J(X_\a)$ with eigenvalue $\mu_\a$. Thus we get
an orthonormal basis $\{ T, X_\a, Y_\a ,\a=1,\dots,n\}$ of $T_PM$.
We may assume $\la_\a\ge \mu_\a$.

The values of $\la_\a$ and $\mu_\a$ are determined by
\eqref{gotto} and by Levi umbilicality. Indeed, letting
$Z_\a=X_\a-i Y_\a$, we have $g_{\a\bb}= 2 \d_{\a\b}$. Since $M$ is
Levi umbilical,
\begin{equation}
\label{pinte}
      2 H = Hg( Z_\a,Z_\aa)
          = g( L(Z_\a),Z_\aa)
         = g( \la_\a X_\a- i\mu_\a Y_\a, X_\a +i Y_\a)
          = \la_\a + \mu_\a.
\end{equation}
Moreover, since
 $h_{\a\b} = g(Z_\a, L(Z_{\b})) =
(\la_\a-\mu_\a) \d_{\a\b}$, it is
$
     h_\a^{\bb} = \frac 12 (\la_\a-\mu_\a) \d_{\a}^{\b}.
$
Thus,
\begin{equation}
 \label{cillo}
    (H^2 -h_{00}H) \d_\a^{\ga}= h_{\a}^{\bb} h_{\bb}^\ga
           =\frac 1 4    (\la_\a-\mu_\a)^2 \d_{\a}^{\la}.
\end{equation}
The solutions to equations \eqref{pinte}  and \eqref{cillo} are $
\la_\a= H + \sqrt{H^2- h_{00} H}$ and $\mu_\a = H - \sqrt{H^2-
h_{00} H}$. The proof of \it Step 3 \rm   is concluded.

\medskip

In \it Step 3\rm, we established that the principal curvatures of
$M$ are the constant numbers \eqref{acquaragia}. By a classical
result  going back to Segre \cite{S}, if a connected hypersurface
in $\R^{N+1}$ has constant principal curvatures, then it must be a
plane, a sphere or a cylinder, i.e.~a Cartesian product $\mathbb
S^p\times\R^{N-p}$, where $\mathbb S^p$ is a  $p$-dimensional
sphere and $0\le p\le N$.
  In particular, a surface with constant
curvatures can have at most two different ones. The numbers in
\eqref{acquaragia} are not pairwise different only in the
following two cases:

Case A: $h_{00} = H$ and

Case B: $h_{00} = 0$.

In Case A, all the principal curvatures are equal to $H$ and the
surface $M$ must be contained in a sphere of radius $\frac 1H$. In
Case B the surface must be  a  cylinder. In the latter case,
equations \eqref{acquaragia} become
\begin{equation} \label{destro}
L(X_\a)=2HX_\a ,\quad  L(Y_\a)=0\quad  \text {and}\quad  L(T)=0.
\end{equation}
Fix a point $P$. After a  complex rotation, we may assume that the
vectors at $P$ satisfying \eqref{destro} are $X_\a = \p_{x_\a}$,
$Y_\a=\p_{y_\a}$ and $T = \p_{y_{n+1}}$. This means that $\ker (L)
= \mathrm{span}\{ \p_{y_h}:h=1,\dots, n+1 \}$. For a cylinder,
$\ker(L)$ is the same at any point (after the trivial
identification between different tangent spaces of $\R^{2n}$).
Moreover, the remaining $n$ principal curvatures are all equal to
$2H$. Then the surface is contained in a cylinder of equation
\[
 \sum_{k=1}^{n+1}(x_k-b_k)^{2} = \frac{1}{4H^2},
 \]
for suitable constants $b_k$.   The proof is concluded.

\end{proof}

\bigskip

\begin{proof} [Proof of Theorem \ref{klino}]
Without loss of generality we can assume $H>0.$

\smallskip
\noindent \it Step A. \rm  First we prove that $h_{\a 0}=0$. Since
$h_{\a\b}=0$,
  \eqref{azzo} becomes
$
 \nabla_\b h_{\a\barga} +
ih_{\a\barga}h_{\b 0}    +2ih_{\a 0}  h_{\b\gg} = 0.
$
Contracting with $g^{\gg\a}$ gives $
 \nabla_\b h_\a^\a  +
i h_\a^\a h_{\b 0}   +2i h_{\a 0}  h_{\b}^\a  = 0. $ The Levi
curvature is constant and then
\begin{equation}
  \label{cod}
   n     H h_{\b0} +2h_{\a 0} h_{\b}^\a=0.
\end{equation}

Denote by $k_{(\la)}$, $\la=1,\dots, n$, the principal Levi
curvatures of $M$ at a point $P$. This means that there is an
orthonormal family of holomorphic vectors $V_{(\la)}=V_{(\la)}^\b
Z_\b \in \H _P$, $\la=1,\dots,n$, such that $h_\b^\a V_{(\la)}^\b
= k_{(\la)} V_{(\la)}^\a$. Contracting \eqref{cod} with
$V_{(\la)}^\b$ yields
\[
     (n H   + 2k_{(\la)}) h(T, V_{(\la)})=0.
\]
By  pseudoconvexity, it is $k_{(\la)}\ge 0$ for all $\la=1,\dots,
n$.  Since $H>0$, this implies $h(T, V_{(\la)})=0$ for any
$\la=1,\dots, n$, which ensures $h_{\a 0}=0$.

Inserting $h_{\a\b}=0$, $h_{\a 0}=0$ and $h_\a^\a = n
H=\,$constant in equations \eqref{azzo}, \eqref{gotta} and
\eqref{noni}, we find
 \begin{subequations}
\begin{eqnarray}
 &
\D_\b h_{\a\gg}=0 \label{pappa},
\\&
  \D_0 h_{\a\bb} = ih_{\a\ll}h_\bb^\ll - ih_{\a\bb}h_{00},
  \label{pippo}
\\&
  \D_\a h_{00}=0.\label{peppe}
  \end{eqnarray}
\end{subequations}

Equation \eqref{peppe}  and $\D T=0 $ imply that $Z h_{00}=0$ for
all holomorphic $Z$. Since the horizontal distribution is bracket
generating,
 we conclude   that $h_{00}$
is constant on $M$. Contracting $\a$ and $\bb$ in \eqref{pippo}
and using $H=\,$constant, we find
$
  h_{\a\ll} h^{\a \ll} = n H h_{00}.
$
If $h_{00}=0$, it follows that $h_{\a\ll}=0$ and thus $H=0$. This
is not possible and $h_{00}$ must be a non zero constant. Since
$h_{\a 0}=0$, by \it Step A \rm we also have $L(T) = h_{00 } T$.

\smallskip\noindent\it Step B. \rm
If $X\in T_P M$ is a real   tangent vector orthogonal to $T$ and
such that $L(X) = \la X$,  then the vector $Y=J(X)$ satisfies
$L(Y) = \mu Y$. This follows from $h_{\a\b}=0$ and can be proved
as in \it Step 2 \rm of the proof of Theorem \ref{teor}. Moreover,
letting $Z= X-i Y$ we have
\[
   0=h(Z,Z) = g( L (Z), Z ) = g(\la X-i\mu Y, X-iY) = \la-\mu.
\]
Therefore $\la=\mu$.

Iterating this process $n$ times, we find an orthonormal basis
 $\{X_\a, Y_\a=J(X_\a), T:  \a=1,\dots,
n\}$ of $T_PM$ such that
\begin{equation}
\label{dario}
  L(X_\a)=\la_\a X_\a,\quad
  L(Y_\a)=\la_\a Y_\a.
\end{equation}
Notice that $L$ sends $\H$ into   $\H$,  because $h_{\a\b}=0$ and
$h_{\a0}=0$. Moreover, letting $Z_{\a} = X_\a - i Y_\a$ we have
$L(Z_\a) = \la_\a Z_\a$. The numbers $\la_1,\dots, \la_n$ are the
eigenvalues of the Levi form at the point $P$, i.e.~$h(Z_\a,
Z_\bb) = \la_\a g(Z_\a, Z_\bb)$.

\smallskip\noindent\it Step C. \rm
We claim that the eigenvalues of $L$ are constant. First observe
that  any pair of points in $M$ can be connected by a horizontal
path $\gamma: [0,1] \to M$, i.e.~a piecewise $C^1$ curve  such
that $g(\dot\ga, T)=0$. This follows from the rank condition
\eqref{rnk}.
 Take $P,Q\in M$ and connect them by a horizontal
curve $\gamma$ with $\gamma(0)=P$ and $\gamma(1)=Q$. Let $\{
X_\a^P, Y_\a^P= J(X_\a^P),T:\a=1,...,n\}$ be an orthonormal basis
of $T_PM$ satisfying \eqref{dario}. Let $Z_\a^P = X_\a^P - i
Y_\a^P$ and let $Z_\a$ be the parallel extension of $Z_\a^P$ along
$\gamma$, that is
\begin{equation}\label{martinez}
 \D_{\dot\g} Z_\a =0\text{ along $\gamma$ and  }Z_\a (P)=
 Z_\a^P.
 \end{equation}
The vector field $Z_\a$ is holomorphic

Equation \eqref{pappa} and its conjugate imply
\begin{equation}
   \label{poppo}
   \D_{\dot \ga}h(Z ,\WW)= 0,
           \quad \text{for all holomorphic $Z, W$.}
 \end{equation}
Then, from  \eqref{poppo} and \eqref{martinez} it follows that
 \[
     \frac{d}{dt}h( Z_\a, Z_\bb ) =
           \D_{\dot\g}h( Z_\a, Z_\bb )
        +h(\D_{\dot\g} Z_\a, Z_\bb)
        +h(  Z_\a,\D_{\dot\g} Z_\bb)=0.
\]
Thus $h( Z_\a, Z_\bb )$ is constant along $\ga$ and
\[
     h( Z_\a, Z_\bb )= h(Z_\a^P, Z_\bb^P) = 2\la_\a\delta_{\a\b},
\]
where the $\la_\a$'s are the Levi eigenvalues at $P$. Since $g$ is
parallel, we also have $g(Z_\a, Z_\bb) = g(Z_\a^P, Z_\bb^P) =
2\d_{\a\b}$. Eventually, we get $h( Z_\a, Z_\bb ) = \la_\a g(Z_\a,
Z_\bb)$,  where the $\la_\a$'s are again the eigenvalues at $P$.
This means that also at the point $Q=\g(1)$ the eigenvalues of $L$
are $\la_1, \la_2, \dots, \la_n$ and $h_{00}$.

\smallskip\noindent\it Step D. \rm
The shape operator $L$ has constant eigenvalues
$\la_1,...,\la_n,h_{00}$. Each eigenvalue $\la_\a$ has
multiplicity 2 and the corresponding eigenspace is a complex
subspace of $\C^{n+1}$. By Segre's theorem on hypersurfaces with
constant curvatures,  $M$ can have no more than two different
constant curvatures and it is contained either in a sphere or in a
cylinder with spherical section. We may assume
$\la_1=\cdots=\la_{m}=0$ and $\la_{m+1}=\cdots=\la_n=h_{00}$ for
some $0\leq m\leq n-1$. In case $m=0$ we have a sphere. In case
$1\leq m\leq n-1$ we have a cylinder of the form \eqref{cyllo}.
The case case $m=n$ is excluded, because we have a cylinder of the
form $\C^n\times\mathbb S^1$ which has $H=0$.

The proof is concluded.
\end{proof}

\end{document}